\documentclass{amsart}

\numberwithin{equation}{section}

\usepackage{xcolor}
\definecolor{darkblue}{rgb}{0.0, 0.0, 0.55}

\usepackage{multicol}
\usepackage{verbatim}
\usepackage{wrapfig}
\usepackage{tcolorbox}
\usepackage{float}

\usepackage[margin=0.6in,paperwidth=8.5in,paperheight=11in]{geometry}

\usepackage{tikz} 
\usetikzlibrary{decorations,calc,arrows,through,intersections,shapes,arrows.meta}

\usepackage[
   colorlinks=true,
  linkcolor=darkblue,  
 citecolor=darkblue,  
   urlcolor=darkblue    
]{hyperref}

\usepackage[figure]{hypcap}



\newtheorem*{theorem*}{Theorem}
\newtheorem*{remark*}{Remark}

\newcommand{\step}[1]{\noindent \emph{Proof of #1. }}

\def\AngleA{30}
\def\AngleB{18}

\begin{document}

\title[Archimedes' Angle Measure]{Archimedes' Angle Measure \\ and \\ { The Fundamental Theorem of Trigonometry }}

\author{Eric Grinberg}
\address{University of Massachusetts, Boston}
\email{eric.grinberg@umb.edu}

\author{Omar Hijab}
\address{Temple University}
\email{omar.hijab@temple.edu}
\date{}

\begin{abstract}
We aim to understand a distant conversation about a trigonometric principle. Although foundational to calculus, this principle and the basic concept of angle measure are often treated tacitly in both elementary and advanced accounts. The primordial idea of counting rotations makes one aware of the integers. In the same vein, Archimedes' measure of partial rotations  hints at the completeness of the real numbers.  By a less well-known implication, these ideas also lead to the complex numbers and  to the Fundamental Theorem of Trigonometry. What does it all mean? And is it true?
\end{abstract}

\maketitle

\begin{multicols}{2}

\section{Introduction}
\subsection{A caf\'e conversation}
\begin{wrapfigure}[8]{r}{0.43\columnwidth}
\vspace{-10pt}
\includegraphics[width=\linewidth]{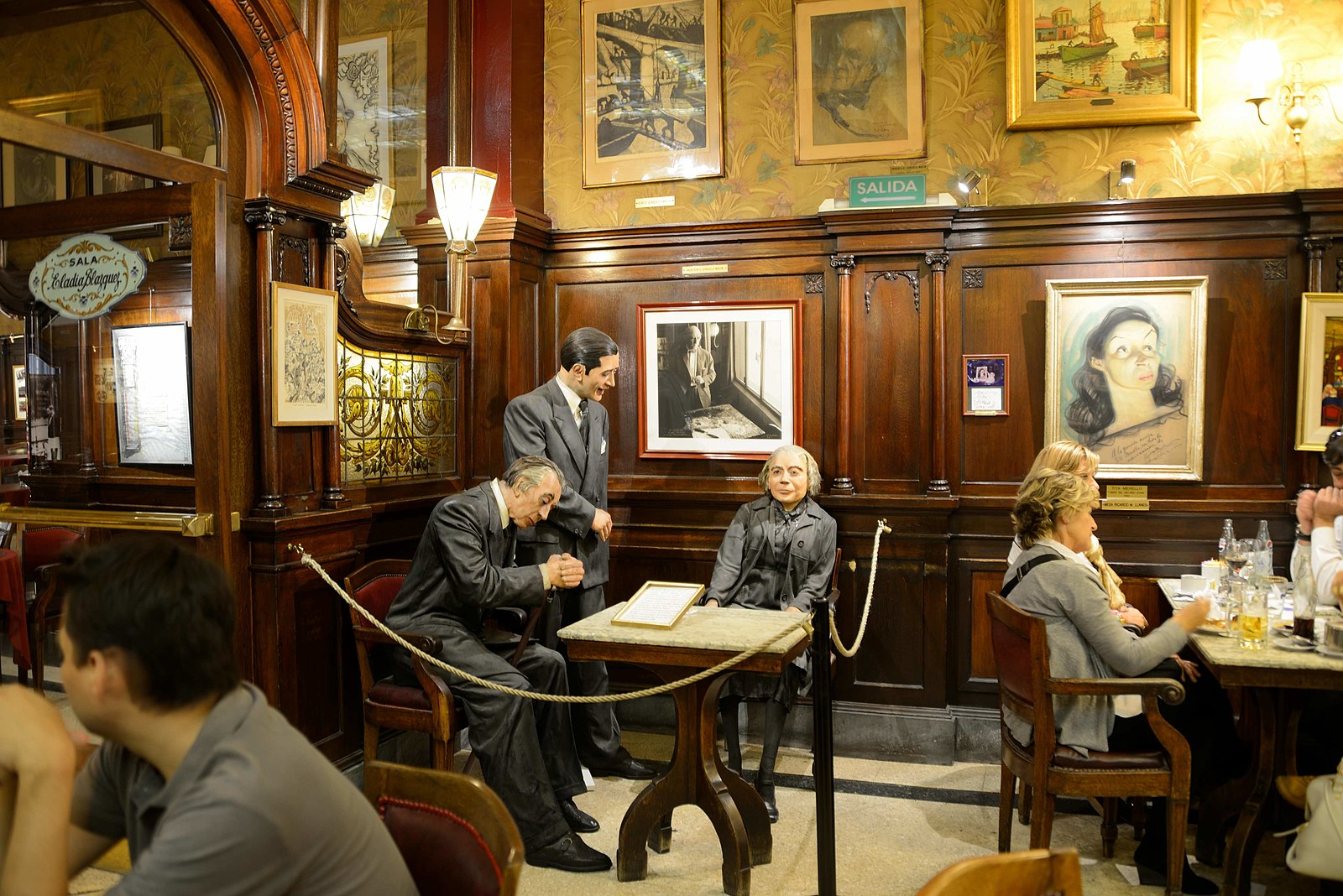}
\end{wrapfigure}
The place was noisy and memories are fuzzy. Maybe the details were different.
Maybe the conversation never took place at all... No matter, we'll try to
recount what we remember and interpret what we recount.

\vspace*{-10pt}
\begin{center} --------\end{center}

S: ``What is the Fundamental Theorem of Trigonometry?"

S: ``Is it
\(
\sin^2\theta  + \cos^2\theta =1 \, ?
\)"

T: ``That's famous and precedes the Fundamental Theorem, but it's a corollary."

S: ``How about
\(
e^{i \theta} = \cos\theta + i \sin\theta \, ? 
\)"

T: ``You're getting warmer, though the coffee is getting colder. This too is a corollary, 
less distant in time from the Fundamental Theorem." 

A server comes and freshens the interlocutors' coffee.

T: (Pauses actively) ``You know,  \emph{\scalebox{1.1}{counting rotations}} is a primordial idea. 
(Pauses again.) It makes us aware of the \emph{\scalebox{1.1}{integers}}. By the same token,
Archimedes' measure of partial rotations leads  to  the \emph{completeness}
of the \emph{\scalebox{1.1}{real numbers}}."  

T: `` $\ldots$ eh, $\ldots$, though less well-known, Archimedes' angle measure also unearths 
the notion of \emph{complex numbers}, 
paving the way to the Fundamental Theorem of Trigonometry."

(T pauses to sip some coffee.)

 T: ``Okay, here it is."

\begin{theorem*}[Fundamental Theorem of Trigonometry]
There is an additive non-constant continuous map  $P(\theta)$  of the real line to the unit circle, unique up to rescaling. This map is surjective.
\end{theorem*}
By \emph{additivity} we mean addition of real numbers corresponds to addition of angles, or equivalently to multiplication on the unit circle \eqref{eq:trighomo}.

\subsection{Stepping back,  seeking perspective}\label{subsec:stepping_back}
Early in the development of trigonometric calculus we come across the identity
\begin{equation}
	\label{eq:sintheta}
	\lim_{\theta\to0}\frac{\sin\theta}\theta=1.
\end{equation} 
But what exactly is $\theta$?

The limit \eqref{eq:sintheta} and its cosine companion are central in developing the derivatives of trigonometric functions. Yet a number of standard textbooks give circular arguments in justifying \eqref{eq:sintheta}, as \cite{Richman} observes.
A mystifying feature of today's calculus books is their lack of discussion of what exactly is meant by the measure $\theta$ of an angle. The subject is assumed to have been treated somehow in secondary school. However, this is rarely so. Secondary school curricula typically follow the G.D.~Birkhoff axioms \cite{B}, and in particular the Birkhoff Protractor Postulate, which assumes the existence of angle measure, bypassing this issue by design. Major texts such as \cite{Ha} and \cite{Lee-2013} are based on the synthetic axiomatization of Hilbert and on a modified set of Birkhoff axioms, respectively.

The notion of \emph{solid angle} may get an entire chapter in an advanced monograph \cite{B-R}. Yet in  pre-calculus or early portions of calculus treatments, instead of developing angle measure  within the framework of just-learned Cartesian geometry,  books often go outside this framework by appealing to pictures not only for motivation, but also for justification. As a consequence, one forfeits the opportunity to present angle measure as a  basic paradigm of calculus.

Typically, calculus books state $\theta$ is the length of the subtended arc along the unit circle, use this to define $\sin\theta$, then go on to derive \eqref{eq:sintheta} by sandwiching a sector's area  between the areas of inscribed and circumscribed triangles. At the beginning of Cartesian geometry, none of this has yet been defined analytically \cite{Richman}.

The book \cite{T} comes close to defining $\theta$, as it discusses the history of the subject at length. The books \cite{Hi}, \cite{L} ignore the issue completely, although the material here would fit perfectly there. 

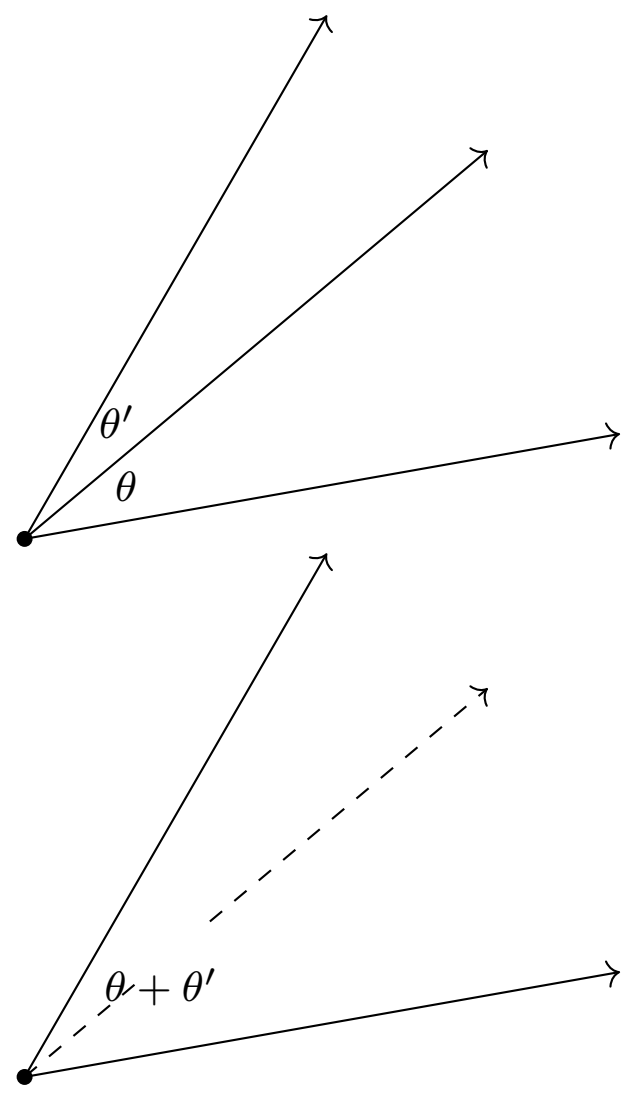
\begin{figure}[H]
\centering
\footnotesize
\begin{tikzpicture}[scale=1.0]
\coordinate (O) at (0,0);
\coordinate (Q) at (10:4);
\coordinate (R) at (40:4);
\coordinate (P) at (60:4);
\draw[->] (O)node[above=10pt, right=14pt]{$\theta$}  -- (Q);
\draw[->] (O)node[above=22pt, right=11pt]{$\theta'$}  -- (R);
\draw[->] (O)  -- (P);
\fill (O) circle [radius=1.5pt];
\end{tikzpicture}

\begin{tikzpicture}[scale=1.0]
\coordinate (O) at (0,0);
\coordinate (Q) at (10:4);
\coordinate (R) at (40:4);
\coordinate (P) at (60:4);
\draw[->] (O)  -- (Q);
\draw[dashed] (O)  -- (40:1);
\draw[->,dashed] (40:1.6)  -- (R);
\draw[->] (O)node[above=17pt, right=12pt]{$\theta+\theta'$}  -- (P);
\fill (O) circle [radius=1.5pt];
\end{tikzpicture}
\caption{Additivity.}
\label{fig:additivity}
\end{figure}

Whatever definition one takes for the measure of the angle, it should be \emph{additive}: When angles are stacked, their measures should add (Figure \ref{fig:additivity}).  

Additive angle measure $\theta$ was used implicitly by Archimedes (287 -- 212 BCE), as part of his work \cite{H} leading to his estimate of the half-circumference of the unit circle,
\begin{equation}
	\label{eq:archest}
	\frac{223}{71}<\pi<\frac{22}7.
\end{equation}
By contrast, Hipparchus (190 -- 120 BCE) \cite{Too1} and Ptolemy (100 -- 170 CE) \cite{Too2} used chord measure $\theta_1$, which is equivalent to $\theta$ but not additive, to build their trigonometric tables. Aryabhata (476 -- 550 CE) presented the  first additive trigonometric tables \cite{SS}. His use of \emph{ardha-jya} for half-chords, shortened to \emph{jya} or \emph{jiva}, morphed into the Arabic \emph{jaib}  (``pocket'') which translates into the Latin \emph{sinus} \cite{MB}.

The path below is in three stages. We first introduce angle stacking and show how it leads naturally toward the complex-number multiplication law. We then follow Archimedes in constructing additive angle measure and proving its existence and uniqueness. Finally, we recover trigonometry from angle measure and establish the Fundamental Theorem of Trigonometry.

\section{Geometry---angles, and how to stack them}
\label{sec:stack}

Working in the Cartesian plane, we will take points to be ordered pairs of numbers, say $P=(x,y)$, $P'=(x',y')$, $P''=(x'',y'')$.   Points may be added, $P + P' := (x+x',y+y')$, and scaled  by numbers $t$, $tP := (tx,ty)$.

An  \emph{angle} is a pair of rays emanating from a point, the vertex of the angle. If the vertex is the origin $O=(0,0)$, then an angle is determined by the intersections of its rays with the unit circle.  If the points of intersection are denoted $N$ and $P$, angles are well-defined even when $P=\pm N$. Notice that angles so-defined are \emph{unoriented}---if we swap rays we obtain the same angle. If we choose an ordering for the rays then we obtain the notion of \emph{oriented angle}. We will be somewhat cavalier about the oriented/unoriented  dichotomy of angles, except where precision is required.

In what follows, we make the normalizing assumption that $N=(1,0)$. Readers needn't worry---generality will not be lost. With this normalization, an angle is determined by a single point $P$ on the unit circle, together with ordering. We need to choose an ordering for the elements of $\{N,P\}$, or eschew orientation.

In the Cartesian plane, angles may be stacked, leading to complex multiplication and division. When stacking two unoriented angles, it's not clear whether we are adding or subtracting them since, for $\angle\alpha$ unoriented, $\angle \alpha = - \angle \alpha$, so $\angle\alpha+\angle\beta$ actually stands for $\pm \angle\alpha \pm \angle\beta$. This is already manifested in ancient tables, such as Ptolemy's \cite{Too2}, where both sums and differences of angles are tabulated.

\begin{figure}[H]
\centering
\resizebox{\columnwidth}{!}{%
    \begin{tikzpicture}[scale=1]
    \coordinate (O) at (0,0);
    \coordinate (P) at (\AngleA:3);
    \coordinate (I) at (3,0);
    \fill[nearly transparent] (O) -- (O-|P) -- (P) -- cycle;
    \draw (P)node[right]{$P$} -- (O-|P);
    \draw (I) -- (O)node[below]{$O$} --node[above]{$1$} (P);
    \draw[fill=black] (P) circle [radius=1pt];
    \draw[fill=black] (O) circle [radius=1pt];
    \end{tikzpicture}%
    \quad
    \begin{tikzpicture}[scale=1]
    \coordinate (O) at (0,0);
    \coordinate (P') at (\AngleB:3);
    \coordinate (I) at (3,0);
    \coordinate (Q) at (O-|P');
    \fill[very nearly transparent] (O) -- (O-|P') -- (P') -- cycle;
    \draw (P')node[right]{$P'$} -- (O-|P');
    \draw (I) -- (O)node[below]{$O$} --node[above]{$1$} (P');
    \draw[fill=black] (P') circle [radius=1pt];
    \draw[fill=black] (O) circle [radius=1pt];
    \draw[fill=white] (Q) circle [radius=1pt];
    \end{tikzpicture}%
}

\resizebox{\columnwidth}{!}{%
    \begin{tikzpicture}[scale=1]
    \coordinate (O) at (0,0);
    \coordinate (P) at (\AngleA:3);
    \coordinate (I) at (0:3);
    \fill[nearly transparent] (O) -- (O-|P) -- (P) -- cycle;
    \draw (I) -- (O)node[below]{$O$} -- (P)node[above right]{$P$} -- (O-|P); 
    \begin{scope}[rotate around={\AngleA:(O)}]
    \coordinate (P') at (\AngleB:3);
    \draw (O) -- (P');
    \coordinate (Q) at (O-|P');
    \fill[very nearly transparent] (O) -- (O-|P') -- (P') -- cycle;
    \draw (Q) -- (P'); 
    \draw[fill=white] (Q) circle [radius=1pt];
    \node[below] at (Q) {$Q\,\,\,$}; 
    \end{scope}
    \draw[fill=black] (O) circle [radius=1pt];
    \draw[fill=black] (P) circle [radius=1pt]; 
    \end{tikzpicture}%
    \quad
    \begin{tikzpicture}[scale=1]
    \coordinate (O) at (0,0);
    \coordinate (P) at (\AngleA:3);
    \coordinate (I) at (3,0);
    \fill[nearly transparent] (O) -- (O-|P) -- (P) -- cycle;
    \draw (P) -- (O-|P);
    \draw (I) arc (0:90:3);
    \draw (I) -- (O)node[below]{$O$} -- (P)node[above right]{$P$}; 
    \begin{scope}[rotate around={\AngleA:(O)}]
    \coordinate (P') at (\AngleB:3);
    \draw (P')node[above]{$P''$};
    \coordinate (Q) at (O-|P');
    \node at (Q) [draw,circle through=(P')] {};
    \draw[fill=black] (P') circle [radius=1pt];
    \draw[fill=white] (Q) circle [radius=1pt] node[below]{$Q$};
    \end{scope}
    \fill[very nearly transparent] (O-|P') -- (P') -- (O);
    \draw[dashed] (O-|P') -- (P') -- (O);
    \coordinate (P'') at (\AngleA-\AngleB:3);
    \draw (P'')node[below=1pt, right]{$P''$};
    \fill[very nearly transparent] (O-|P'') -- (P'') -- (O);
    \draw[dashed] (O-|P'') -- (P'') -- (O);
    \draw[fill=black] (P'') circle [radius=1pt];
    \draw[fill=black] (O) circle [radius=1pt];
    \draw[fill=black] (P) circle [radius=1pt];
    \end{tikzpicture}%
}

\resizebox{\columnwidth}{!}{%
    \begin{tikzpicture}[scale=1]
    \coordinate (O) at (0,0);
    \coordinate (P) at (\AngleA+\AngleB:3);
    \coordinate (I) at (3,0);
    \fill[nearly transparent] (O) -- (O-|P) -- (P) -- cycle;
    \draw (P)node[right]{$P''$} -- (O-|P);
    \draw (I) -- (O)node[below]{$O$} --node[above]{$1$} (P);
    \draw[fill=black] (P) circle [radius=1pt];
    \draw[fill=black] (O) circle [radius=1pt];
    \end{tikzpicture}%
    \quad
    \begin{tikzpicture}[scale=1]
    \coordinate (O) at (0,0);
    \coordinate (P'') at (\AngleA-\AngleB:3);
    \coordinate (I) at (3,0);
    \coordinate (Q) at (O-|P');
    \fill[nearly transparent] (O) -- (O-|P'') -- (P'') -- cycle;
    \draw (P'')node[right]{$P''$} -- (O-|P'');
    \draw (I) -- (O)node[below]{$O$} --node[above]{$1$} (P'');
    \draw[fill=black] (P'') circle [radius=1pt];
    \draw[fill=black] (O) circle [radius=1pt];
    \end{tikzpicture}%
}
\caption{Angle stacking.}
\label{fig:stacking}
\end{figure}
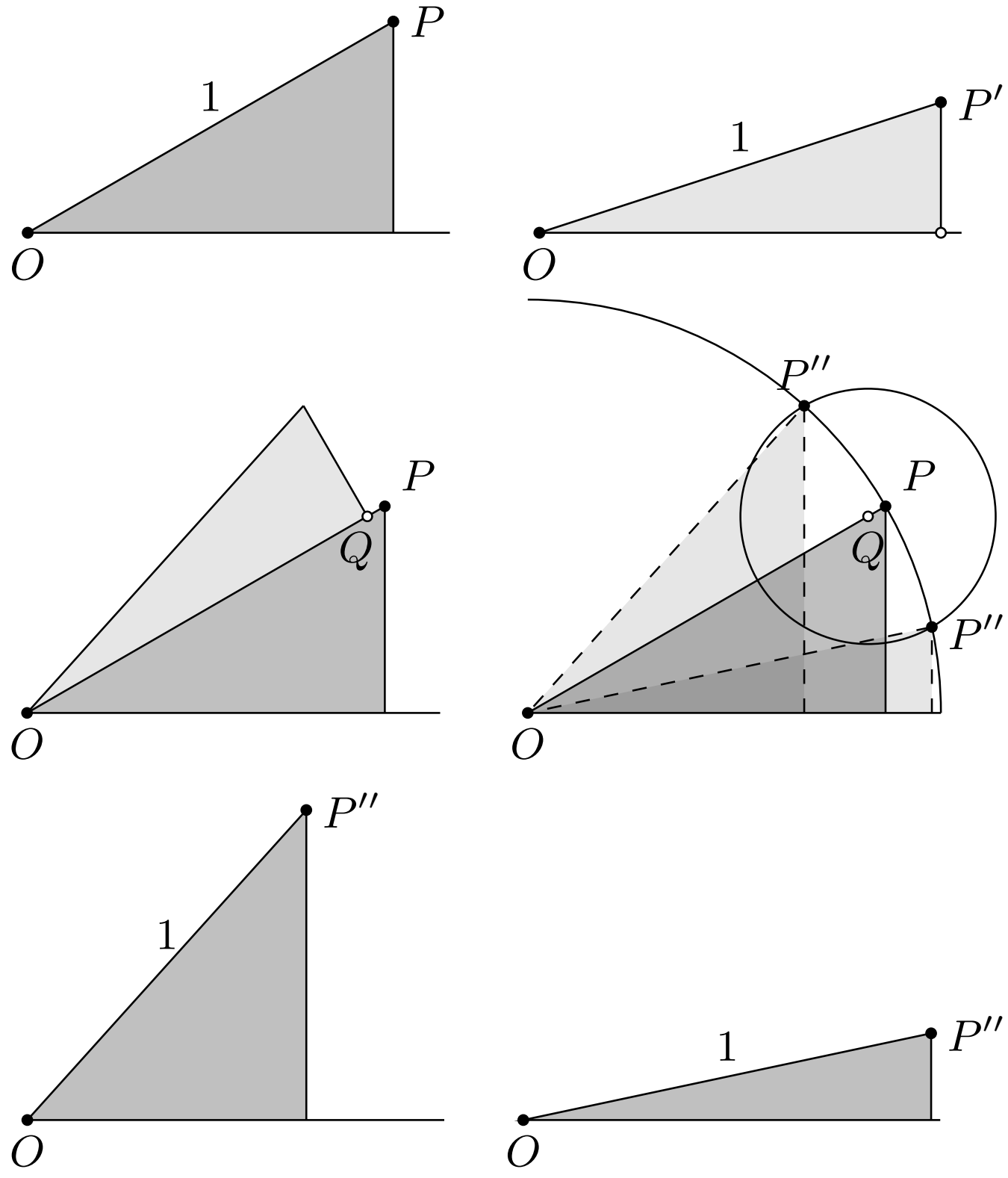

Let $P$ and $P'$ be points on the unit circle corresponding to two given angles. Then the angles may be stacked (Figure \ref{fig:stacking}), resulting in the angle corresponding to a point $P''$ on the unit circle.  The illustration seems to suggest subliminal orientation, denoting a sum and not a difference. But not to worry, the geometry and the algebra will ally to bring out both possibilities together.

To make this precise, draw the circle with center $Q=x'P$ and radius $|y'|$. Note that in the figure, $x'$ is only slightly smaller than $1$, the length of the hypotenuse of triangle $P$. The small gap is significant, as $Q$ and $P$ play distinct roles.

When $x'y' \neq 0$, we claim this circle intersects the unit circle at two points, both labelled $P''$, for {\em beneficial ambiguity}:
\footnote{This phrase is in good company with \emph{abstract nonsense} and \emph{abuse of notation}. See also the definition of \emph{beneficial insects} in the Old Farmer's Almanac; this is only jargon, not meant to denigrate any insects.} 

\enlargethispage{\baselineskip}
\begin{equation} 
\label{eq:product}
\begin{aligned}
P''=&(xx'- yy', x'y+xy') \\ P''=&(xx'+ yy', x'y-xy').
\end{aligned} 
\end{equation}


Since we obtain two points $P''$, angle stacking in fact results in two angles, which may remind the reader of the complex product and complex quotient of $P$ and $P'$. Although we do not develop complex numbers here, we'll speak of the \emph{real} and \emph{imaginary} parts of a point in the plane.

To derive \eqref{eq:product}, let $\langle P,P'\rangle= xx'+yy'$, $|P|^2=\langle P,P\rangle$, and $P^\perp=(-y,x)$. Then $P$ is on the unit circle if and only if $|P|^2=1$, and $P$, $P'$ on the unit circle satisfy  $\langle P,P'\rangle=0$ if and only if $P'=\pm  P^\perp$.

Since $P''$ is on the circle of center $Q$ and radius $|y'|$, we may write $P''=Q+y'R$ for some point $R$ on the unit circle. Then  
\begin{equation*}
\begin{split}
1 =|P''|^2 &=|x'P+y'R|^2 \\
&={x'}^2|P|^2+ {y'}^2|R|^2+2x'y'\langle P,R\rangle\\
&= ({x'}^2+ {y'}^2)+2x'y'\langle P,R\rangle\\
&=1+2x'y'\langle P,R\rangle,
\end{split}
\end{equation*}
which occurs if and only if $\langle P,R\rangle=0$, i.e., if and only if  $R=\pm P^\perp$. This  yields $P''=x'P\pm y'P^\perp$, which is \eqref{eq:product}.   That establishes the claim.

Given any two points $P=(x,y)$, $P'=(x',y')$ on the unit circle, we can define their product $PP'$ and quotient $P/P'$ by \eqref{eq:product}. One can check directly that this multiplication rule is commutative and associative, with $N$ as an identity element: for any point $R$ on the unit circle,
\[
PP'=P'P,\quad P(P'R)=(PP')R,\quad  PN=NP=P.
\]
(Just copy and paste from an introduction to complex numbers, but don't tell.) 
If $\bar P := (x,-y)$, then $P\bar P=N$ for any point $P$ on the unit circle, and $P/P'=P\overline {P'}$ is the Hermitian product.

In the next section,  we define an additive angle measure $\theta(P)$, assigning real numbers to points $P$ in the punctured circle $\{ x\ne -1\}$. To take into account both cases in \eqref{eq:product}, we require 
\begin{equation}
\label{eq:triginv}
\theta(P)=-\theta(\bar P),\quad \text{for} \quad x\ne-1,
\end{equation}
and consequently, $\theta(N)=0$. We then interpret Figure \ref{fig:additivity} by 
\begin{equation}
\label{eq:additivity}
\theta(PP')=\theta(P) + \theta(P'),\qquad x>0,x'>0.
\end{equation}
It is in this form that we shall establish additivity of $\theta=\theta(P)$.

Note that \eqref{eq:additivity} implies
\begin{equation}
\label{eq:trigdouble}
\theta(P^2)=2\theta(P),\qquad x>0.
\end{equation}

For $P=(x,y)$ on the punctured unit circle $x\ne-1$, write
\begin{equation}
\label{eq:mean}
\sqrt{P}= \left(\frac{1+x}{\sqrt{2+2x}},\frac{y}{\sqrt{2+2x}}\right).
\end{equation}
Then $P\mapsto \sqrt{P}$ maps the punctured unit circle continuously and bijectively onto the right-half unit circle $x>0$, with inverse $P\mapsto P^2$, and the imaginary parts of $P$ and $\sqrt{P}$ have the same sign.

Since $(\sqrt{P}\sqrt{P'})^2=PP'$, we have
\begin{equation}
\label{eq:sqrt}
\sqrt{PP'}=\sqrt{P}\sqrt{P'},
\end{equation}
up to sign, whenever both sides of \eqref{eq:sqrt} are defined. (This is one instance where we benefit from beneficial ambiguity.) When the imaginary parts of $P$, $P'$, and $PP'$  have the same sign,   the imaginary parts of the sides of \eqref{eq:sqrt} have the same sign. Hence \eqref{eq:sqrt} is correct as written, when  the imaginary parts of $P$, $P'$, and $PP'$  have the same sign.

\section{Measurement} 
\subsection{Enter Archimedes}
Archimedes interpreted angle measure $\theta$ as the length of the arc along the unit circle subtended by the angle (Figure \ref{fig:archpost}). 
In his treatises \emph{On the sphere and the cylinder} and \emph{Measurement of a circle} \cite{H}, Archimedes grappled with the ideas underlying the length, area, and volume of curves, surfaces, and solids.  

Archimedes not only computed the circumference and area of a circle, the surface area of a sphere, and the volume of a ball, he also isolated  \emph{a priori} properties such measures ought to have. This led to a plethora of ideas including, explicitly or implicitly,  the definition of a \emph{convex} curve, the Archimedes Convexity Postulate, the additivity of angle measure (\S\ref{sec:stack}), the bisection method (\S\ref{sec:bisection}), and the completeness of the real number line. 

The \emph{Archimedes Convexity Postulate} (ACP) \cite{E}, \cite{GW}, states:  
 
 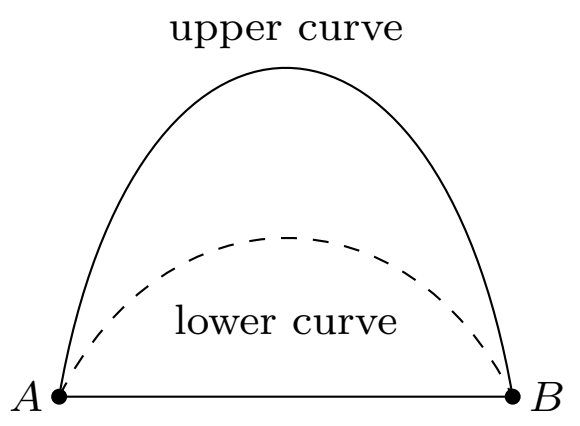
\begin{figure}[H]
\footnotesize
\centering
\begin{tikzpicture}[scale=1.0]
\coordinate (A) at (0,0);
\coordinate (B) at (3,0);

\draw (A) node[left] {$A$} -- (B) node[right] {$B$};

\draw[dashed] (A) .. controls (0.7, 1.4) and (2.3, 1.4) .. (B);
\node[below=8pt] at (1.5, 1.0) {lower curve};

\draw (A) .. controls (0.5, 2.9) and (2.5, 2.9) .. (B);
\node[above] at (1.5, 2.2) {upper curve};

\fill (A) circle [radius=1.5pt];
\fill (B) circle [radius=1.5pt];
\end{tikzpicture}
\caption{Archimedes Convexity Postulate (ACP).}
\label{fig:convex_postulate}
\end{figure}
 
\begin{tcolorbox}[colback=gray!30, colframe=gray!30, arc=2mm, boxrule=0.5pt]
\small\itshape
If two convex curves have common endpoints $A, B$ along a base segment and 
one curve lies below the other then the lower curve has smaller length than the upper curve.
\end{tcolorbox}

Let $P=(x,y)$ be a point on the unit circle first quadrant. In Figure~\ref{fig:archpost}, $\theta$ is the as-yet-undefined length of the arc. We will use the ACP to derive the following estimates:
		\begin{equation}
			\label{eq:archpost}
			y <\theta<1-x+y,\quad x,y>0 \quad 
			\text{\small \emph{(Elbow Bounds)}}.
		\end{equation}
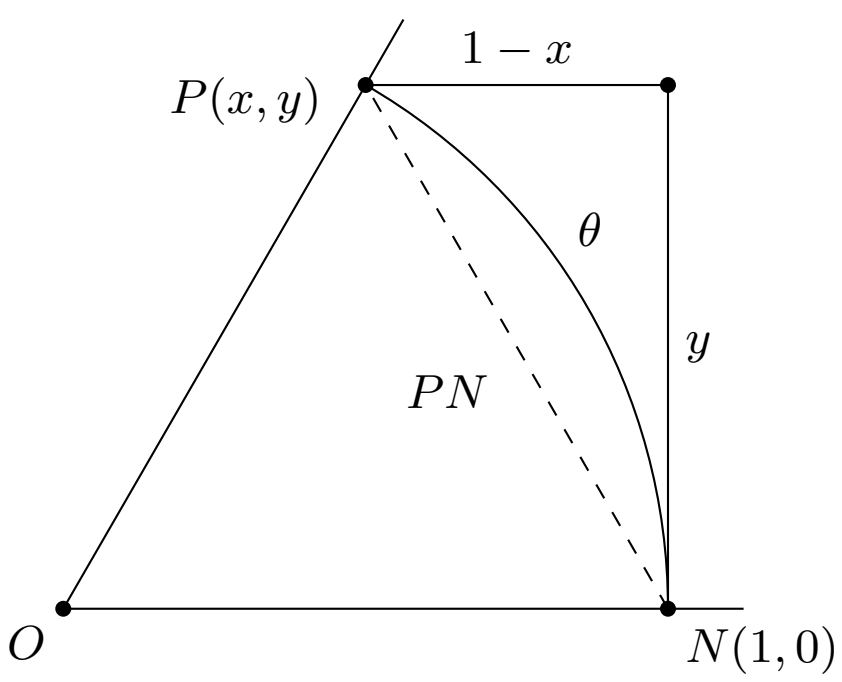
\begin{figure}[H]
\begin{tikzpicture}[scale=1.0]
\coordinate (O) at (0,0);
\coordinate (N) at (0:4);
\coordinate (P) at (60:4);

\draw (O) node[below left]{$O$} -- (4.5,0);
\draw (O) -- (60:4.5);

\node[left=5pt, yshift=-3pt] at (P) {$P(x,y)$};

\draw (N) arc (0:60:4);
\node[above right=-2pt] at (35:4.1) {$\theta$};

\draw[dashed] (N) node[below right]{$N(1,0)$} -- node[midway, below left=2pt] {$PN$} (P);

\draw (N) -- node[midway, right] {$y$} (N |- P) -- node[midway, above] {$1-x$} (P);
\fill (N |- P) circle (1.5pt);

\fill (O) circle [radius=1.5pt];
\fill (N) circle [radius=1.5pt];
\fill (P) circle [radius=1.5pt];
\end{tikzpicture}
\caption{Elbow Bounds.}
\label{fig:archpost}
\end{figure}

In Figure \ref{fig:archpost}, applying the ACP (with base $PN$) to the curve $\theta$ (the lower curve) and the elbow consisting of $1-x$ and $y$ (the upper curve), we obtain
\[
\theta<1-x+y \, ;
\]
applying the postulate to the chord $PN$ and the arc $\theta$ we obtain $PN < \theta$. Then, by the Pythagorean theorem, $y < PN$. Putting these together we have \eqref{eq:archpost}.

This is all we require of the ACP and its direct role ends here: the postulate serves to certify that the Elbow Bounds capture what arc length ought to satisfy. Below we derive the existence and uniqueness of an additive angle measure $\theta$ satisfying the Elbow Bounds.

\subsection{Archimedean bisection: chord estimates}\label{sec:bisection}
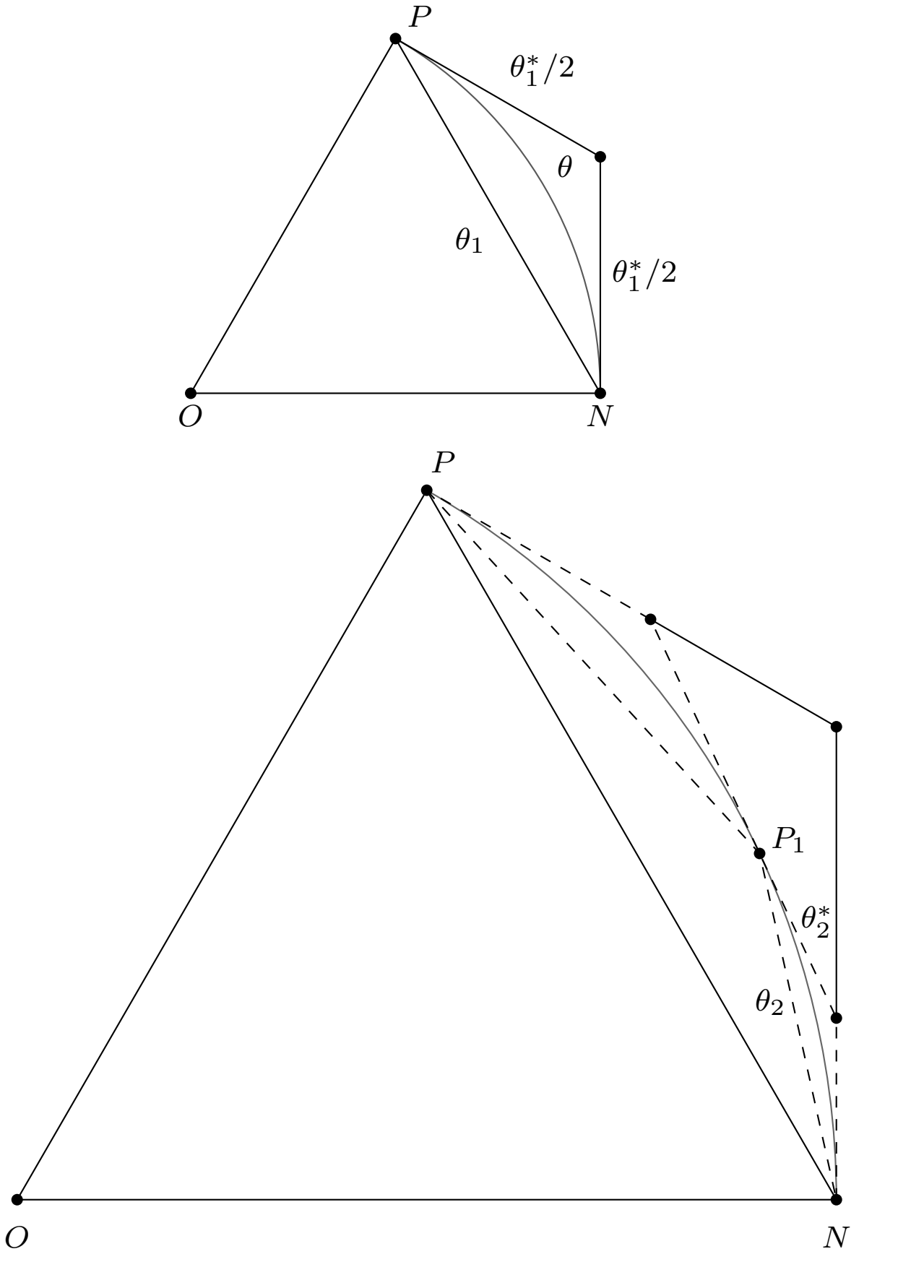
\begin{figure}[H]
			\centering
			\footnotesize
			
			\def\a{60}
			\def\b{25}
			\begin{tikzpicture}[scale=0.75]
				\coordinate (O) at (0,0);
				\coordinate (I) at (0:5);
				\coordinate (P) at (\a:5);
				\coordinate (T1) at (5,{5*(1-cos(\a))/sin(\a)});
				\draw[black!65] (I) arc (0:\a:5);
				\draw   (O)node[below]{$O$} --  (I)node[below]{$N$} -- (T1) --(P) -- cycle;
				\draw (P)node[above right]{$P$} -- (I);
				\fill (P) circle [radius=2pt];
				\fill (I) circle [radius=2pt];
				\fill (O) circle [radius=2pt];
				\fill (T1) circle [radius=2pt];
				\node[above right] at (30:5) {$\theta$};                 
				\path (I) -- (P)  node[midway, below left]{$\theta_1$};  
				\path (T1) -- (P)  node[midway, above right]{$\theta_1^*/2$};  
				\path (I) -- (T1) node[midway, right]{$\theta_1^*/2$};      
			\end{tikzpicture}
			
			\begin{tikzpicture}[scale=1.5]
				\coordinate (O) at (0,0);
				\coordinate (I) at (0:5);
				\coordinate (P) at (\a:5);
				\draw (\a:5.2)node{$P$};
				\coordinate (P1) at (\b:5);
				\draw (\b:5.2)node{$P_1$};
				\coordinate (T) at (\b:{5*sec(\b)});
				\coordinate (T1) at (5,{5*(1-cos(\a))/sin(\a)});
				\coordinate (T2) at (5,{5*(1-cos(\b))/sin(\b)});
				\coordinate (T2p) at (intersection of  P1--T2 and P--T1);
				\draw[black!60] (I) arc (0:\a:5);
				\draw (P) -- (I);
				\draw (O)node[below, yshift=-4pt]{$O$} -- (I)node[below, yshift=-4pt]{$N$};
				\draw (O) -- (P);                               
				\draw (T2) -- (T1) -- (T2p);
				\draw[dashed] (I) -- (T2) -- (T2p) -- (P);
				\draw[dashed] (P) -- (P1) -- (I);
				\fill (P) circle [radius=1pt];
				\fill (I) circle [radius=1pt];
				\fill (P1) circle [radius=1pt] ;
				\fill (O) circle [radius=1pt];
				\fill (T1) circle [radius=1pt];
				\fill (T2) circle [radius=1pt] ;
				\fill (T2p) circle [radius=1pt] ;
				\path (I) -- (P1) node[midway, above left]{$\theta_2$};  
				\node[above = 18pt, xshift=-5pt] at (T2) {$\theta_2^*$};                       
			\end{tikzpicture}
			\caption{Inscribed and circumscribed \\ chords; lower diagram enlarged.}
			\label{fig:inscirc}
		\end{figure}

\begin{figure}[H]
			\centering
			\footnotesize
			\begin{tikzpicture}[scale=1.0]
				\coordinate (O) at (0,0);
				\coordinate (I) at (0:5);
				\coordinate (P) at (70:5);
				\draw (70:5.3)node{$P$};
				\coordinate (P1) at (35:5);
				\draw (35:5.3)node{$P_1$};
				\coordinate (P2) at (17.5:5);
				\draw (17.5:5.3)node{$P_2$};
				\coordinate (M2) at (intersection of  O--P2 and I--P1);
				\draw[black!70] (I) arc (0:70:5);
				\draw (O) -- (I);
				\draw (O)node[below]{$O$} -- (P);
				\draw (P2) --node[above=3pt, right]{ $\frac14\theta_3$} (I);
				\draw (P1) --node[above=8pt, left]{ $\frac12\theta_2$}  (I);
				\draw (P) --node[below=3pt, left]{ $\theta_1$}  (I);
				\draw (I)node[below right]{$N$};
				\fill (P) circle [radius=1.5pt];
				\fill (P1) circle [radius=1.5pt];
				\fill (P2) circle [radius=1.5pt];
				\fill (I) circle [radius=1.5pt];
				\fill (O) circle [radius=1.5pt];
			\end{tikzpicture}
			\caption{Bisection. }
			\label{fig:bisection}
		\end{figure}
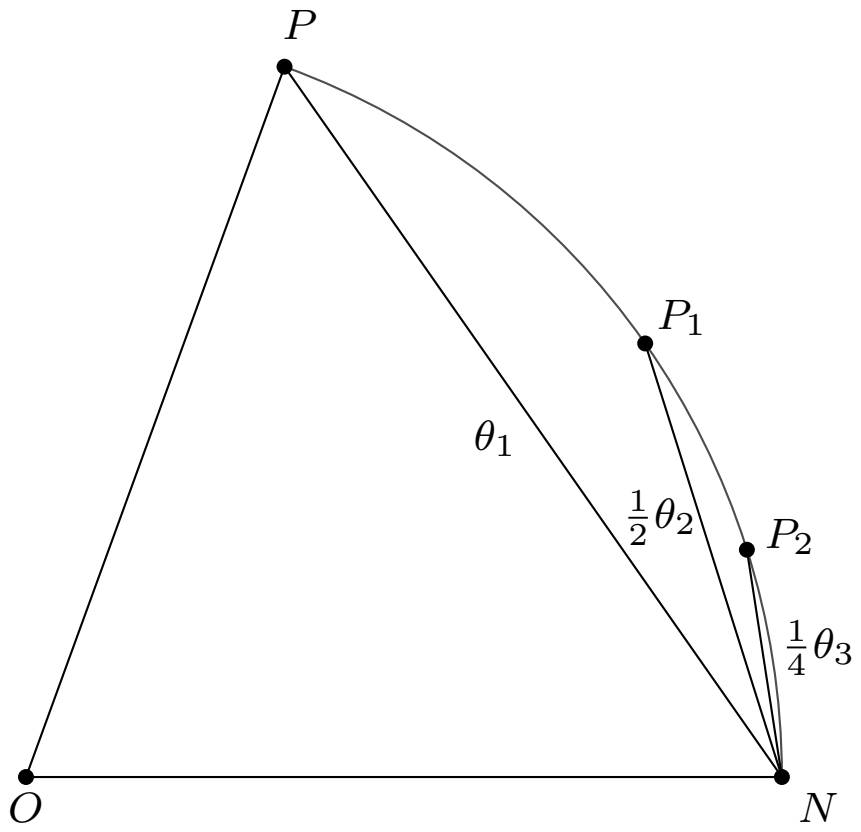

\begin{figure}[H]
			\centering
			\footnotesize
			\def\a{60}
			\def\b{30}
			\def\c{15}
			\begin{tikzpicture}[scale=1.0]
				\coordinate (O) at (0,0);
				\coordinate (I) at (0:5);
				\coordinate (P) at (\a:5);
				\draw (\a:5.3)node{$P_{n-1}$};
				\coordinate (P1) at (\b:5);
				\coordinate (P2) at (\c:5);
				\coordinate (Q) at (O-|P);
				\coordinate (Q1) at (O-|P1);
				\coordinate (T) at (5,{5*sin(\a)}); 
				\draw (T)node[right]{$R_{n-1}$};
				\coordinate (T1) at (\b:{5*sec(\b)});
				\coordinate (M1) at (intersection of  O--P1 and I--P);
				\coordinate (T2) at (\c:{5*sec(\c)});
				\fill (T2) circle [radius=1.5pt];
				\coordinate (M2) at (intersection of  O--P2 and I--P1);
				\fill (M2) circle [radius=1.5pt];
				\fill (P2) circle [radius=1.5pt];
				\draw[dashed] (O) -- (T2)node[right]{$T_{n+1}$};
				\draw (P1) -- (T2);
				\coordinate (Q2) at (O-|P2);
				\draw (O) -- (T2)node[right]{$T_{n+1}$};
				\draw[black!70] (I) arc (0:\a:5);
				\draw (O) -- (M1)node[left=13pt,above=-2pt]{$M_{n}$} -- (T1)node[right]{$T_{n}$};
				\draw (O)node[below]{$O$} -- (P1)node[right=11pt,below=-4pt]{\scriptsize $P_{n}$};
				\draw (I)node[below right]{$N$} -- (P);
				\draw  (O) -- (I) -- (T1) -- (P) -- cycle;
				\draw[dashed]  (P) -- (T) -- (T1);
				\draw (Q1)node[below]{$Q_{n}$};
				\draw (Q)node[below]{$Q_{n-1}$};
				\draw (P) -- (P1) -- (I);
				\fill (P) circle [radius=1.5pt];
				\fill (M1) circle [radius=1.5pt];
				\fill (I) circle [radius=1.5pt];
				\fill (P1) circle [radius=1.5pt] ;
				\fill (O) circle [radius=1.5pt];
				\fill (T) circle [radius=1.5pt];
				\fill (T1) circle [radius=1.5pt];
				\fill (Q) circle [radius=1.5pt];
				\fill (Q1) circle [radius=1.5pt];
				\fill (Q2) circle [radius=1.5pt];
			\end{tikzpicture}
			\caption{First Cut. }
			\label{fig:archdiag}
		\end{figure}
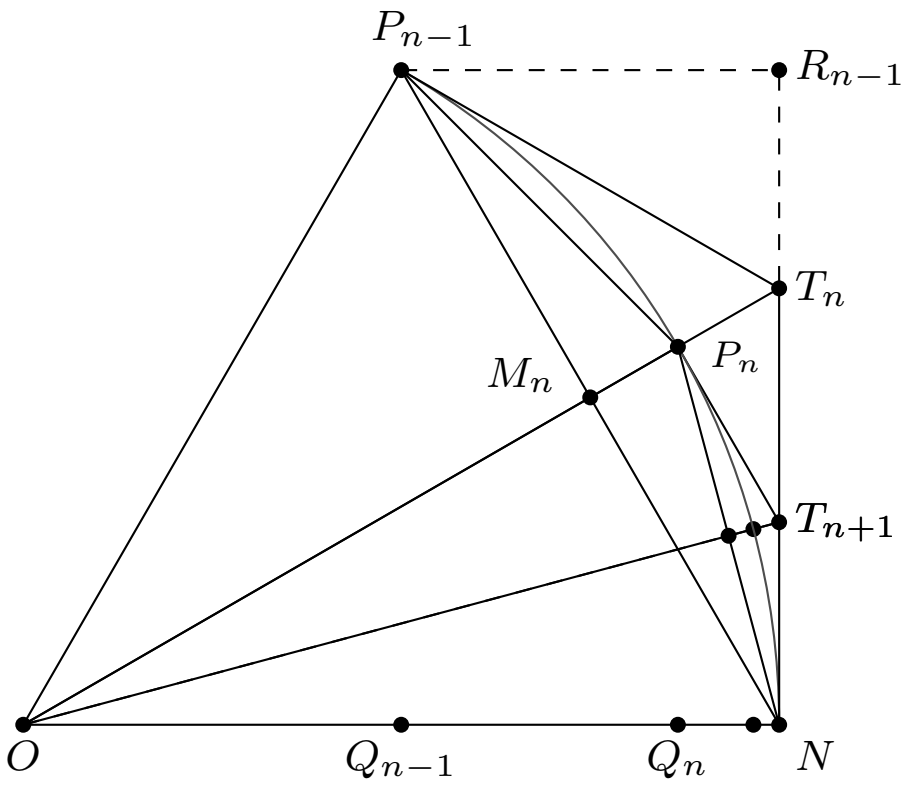

		Archimedean bisection is based on  Figure \ref{fig:inscirc},  a starting point for the \emph{method of exhaustion} of Eudoxus and Euclid. Figures \ref{fig:inscirc}, \ref{fig:bisection}, and \ref{fig:archdiag} are in Archimedes' \emph{On the sphere and the cylinder}  (see the figures of Propositions 1--6 therein) and his \emph{Measurement of a circle}.   
		
		In Figure \ref{fig:inscirc}, let $N$, $P$, $P_1$ be three points on the upper-half unit circle, $y>0$. Let $\theta_1$, $\theta_1^*$ be the inscribed and circumscribed chord-length sums on the top, and let $\theta_1$, $\theta_2$, $\theta_2^*$, $\theta_1^*$
		be the inscribed and circumscribed chord-length sums on the bottom. By the triangle inequality, or by the ACP applied to the four polygonal convex curves joining the endpoints $N$ and $P$, 
		\[
		\theta_1<\theta_2<\theta_2^*< \theta_1^*,\qquad y>0.
		\]
		By incrementally adding points on the arc, partitioning subarcs between adjoining points, we obtain inscribed and circumscribed polygonal approximations of the arc whose lengths are arbitrarily close to each other. This is Proposition 3 in \emph{On the sphere and the cylinder}.
		By adding points repeatedly to the arc, we will establish the Archimedes sequences
		\begin{equation}
			\label{eq:archseq}
			\theta_1<\theta_2<\theta_3<\dots< \theta_3^* < \theta_2^* < \theta_1^*,\qquad y>0.
		\end{equation}
		This is for any point $P$ in the upper-half circle, while the Elbow Bounds are only for $P$ in the unit circle first quadrant. 
		
		To obtain \eqref{eq:archseq} assume the points $P=(x,y)$, $P_1$, $P_2$, $P_3$, $\dots$ are generated by bisection (Figure~\ref{fig:bisection}): $P_1$ bisects the arc from $N$ to
		$P$, and each $P_n$ bisects the arc from $N$ to $P_{n-1}$. Thus $P_n$
		cuts off $1/2^{n}$ of the arc from $N$ to $P$, and the corresponding
		chords $NP_n$ shrink accordingly. In particular, with $P_n=(x_n,y_n)$, the sequences $x_n$ and $y_n$ are increasing and decreasing respectively.
		
		We now focus on the $n$-th bisection. In Figure~\ref{fig:archdiag}, let $M_n$ be the midpoint of the chord $P_{n-1}N$, and let $P_n$ and $T_n$ be the intersections of the line through $O$ and $M_n$ with the unit circle and the vertical line through $N$. Let $P_{n-1}=(x_{n-1},y_{n-1})$, $R_{n-1}=(1,y_{n-1})$, $P_n=(x_n,y_n)$, $Q_{n-1}=(x_{n-1},0)$, $Q_n=(x_n,0)$, $T_n=(1,y_n/x_n)$. Then $P_n$
		is well-defined if  $y_{n-1}>0$, and computing $M_n$ then $P_n=M_n/|M_n|$ shows $P_n=\sqrt{\mathstrut P_{n-1}}$.
		
		From the Figures, it makes sense to define $\theta_n := 2^{n-1}\vert P_{n-1}N \vert$ and $\theta_n^* := 2^n\vert T_nN \vert$. Since $M_n$ is the chord $P_{n-1}N$'s  midpoint,  $\theta_n=2^n\vert M_nN \vert $.  Since the triangles $OM_nN$ and $OP_nQ_n$ are congruent, $\vert M_nN\vert =\vert P_nQ_n\vert =y_n$, hence $\theta_n=2^ny_n$, which implies $\theta^*_n=2^ny_n/x_n=\theta_n/x_n$. Thus, for $n\ge1$ and $y>0$, we conclude
		\begin{equation}
			\label{eq:2^ny}
			P=P_n^{2^n},\quad P_n=\sqrt[2^n]{P},\quad\theta_n^*=\frac{\theta_n}{x_n},\quad \theta_n = 2^ny_n.
		\end{equation}
		%
		Moreover $2^{n}\vert P_{n}N\vert > 2^n \vert M_nN \vert$ and 
        \begin{equation*} 
        \begin{split} 
        2^{n+1}\vert T_{n+1}N\vert  &= 2^n(\vert P_nT_{n+1}\vert +\vert T_{n+1}N \vert ) \\
        &< 2^n(\vert T_nT_{n+1}\vert +\vert T_{n+1}N) \vert = 2^n\vert T_nN\vert ,
        \end{split}
        \end{equation*}
        hence
		$$\theta_n<\theta_{n+1}<\theta_{n+1}^*<\theta^*_{n},$$
		which implies \eqref{eq:archseq}.

\subsection{Uniqueness of angle measure}
{\em There is at most one real-valued additive angle measure  on the punctured unit circle satisfying the Elbow Bounds.}

		To establish uniqueness of angle measure, let $\alpha=\alpha(P)$ and $\beta=\beta(P)$ be real-valued additive angle measures  satisfying the Elbow Bounds \eqref{eq:archpost}. By additivity \eqref{eq:additivity}, it is enough to show $\alpha=\beta$ for  $P$ in the unit circle first quadrant. With $P_n$ as above, let $\alpha_n=\alpha(P_n)$ and $\beta_n=\beta(P_n)$ be the corresponding angle measures. Replacing $P$ by $P_n=(x_n,y_n)$ in \eqref{eq:archpost},
		\[
		y_n  < \alpha_n< 1-x_n+y_n,\qquad
		y_n < \beta_n < 1-x_n+y_n.
		\]
		By additivity \eqref{eq:additivity}, $2^n\alpha_n=\alpha$ and $2^n\beta_n=\beta$.  This implies
		\[
		\left|\alpha-\beta\right| =2^n|\alpha_n-\beta_n| < 2^n(1-x_n) = 2^n\frac{y_n^2}{1+x_n} < 2^{-n}\alpha^2
		\]
		for $n\ge1$. Hence $\alpha=\beta$, establishing uniqueness.

\begin{remark*} Note that we did not require continuity. Rigorous treatments of \eqref{eq:sintheta} assume a smoothness context from the start---via power series \cite{Rudin, Landau-1956}, or a differential equation and arclength integral \cite{Hardy-1908,CJ}; continuity follows as a consequence of smoothness. Here it is derived from additivity and the Elbow Bounds.
\end{remark*}

\subsection{Existence of angle measure}
{\em 	There is at least one  real-valued additive angle measure  on the punctured unit circle satisfying the Elbow Bounds \eqref{eq:archpost}.}

Turning to existence,  the sequences in \eqref{eq:archseq} are bounded. By \eqref{eq:2^ny},  $y_n\to0$, hence $x_n\to1$, as $n\to\infty$.
Since  the real line is complete, the sequences in \eqref{eq:archseq} have a common limit $\theta=\theta(P)$ for $P$ in the upper-half circle $y>0$.
Extend  $\theta(P)$ to the lower-half unit circle $y<0$ by  $\theta(P)=-\theta(1/P)$, and set $\theta(N)=0$.  Then, by construction,  \eqref{eq:triginv} and \eqref{eq:trigdouble}  are valid. From Figure~\ref{fig:archdiag} with $n=1$ and $P_0=P$, 
\begin{equation*} 
\begin{split} 
y < \vert PN\vert  &=\theta_1 < \theta<\theta^*_1 = \vert PT_1\vert +\vert T_1N\vert \\
  &<\vert PR_0\vert +\vert R_0T_1\vert +\vert T_1N\vert =1-x+y,
  \end{split}
  \end{equation*}
verifying the elbow bounds \eqref{eq:archpost}. 
 
 Let $P''=PP'=(x'',y'')$, and let $\theta=\theta(P)$, $\theta'=\theta(P')$, and $\theta''=\theta(P'')$. Let $P_n=(x_n,y_n)$, $P'_{n}=(x'_{n},y'_{n})$, $P''_{n}=(x''_{n},y''_{n})$, $n\ge1$, be the bisection sequences starting from   $P$, $P'$, $P''$ respectively, with  $\theta_n$, $\theta'_{n}$, $\theta''_{n}$, $n\ge1$,  the corresponding chord-length sums. It remains to verify the additivity property \eqref{eq:additivity}.

If $yy'y''=0$, then \eqref{eq:additivity} is immediate, hence we may assume $yy'y''\not=0$. Then there are two cases, $yy'>0$ and $yy'<0$. 

If $yy'>0$, by \eqref{eq:triginv}, \eqref{eq:additivity} is valid for $P$, $P'$ if and only if \eqref{eq:additivity} is valid for $\bar P$, $\bar P'$. Hence we may assume $P$, $P'$ are in the first quadrant. It follows $P_n$, $P'_n$ are in the first quadrant for $n\ge1$.
By induction and \eqref{eq:sqrt},  $P''_{n}=P_{n}P'_{n}$. By \eqref{eq:product}, 
$y''_{n}=x'_{n}y_{n}+x_{n}y'_{n}$. By \eqref{eq:2^ny}, 
\begin{equation}
\label{eq:half-chord_addition}
\theta''_{n}=x'_{n}\theta_{n}+x_{n}\theta'_{n} \qquad \text{(half-chord addition identity)}.
\end{equation}
This identity is the Hipparchus–Ptolemy–Aryabhata chord addition alluded to in \S\ref{subsec:stepping_back}. 
Since   $x_{n} \to 1$, $x'_{n} \to 1$,  we obtain $\theta''=\theta + \theta'$, which is \eqref{eq:additivity}. 

If $yy'<0$, then  $x''=xx'-yy'>0$. Since $yy'<0$, $y''(-y)$ and $y''(-y')$ have opposite signs. By switching the roles of $P$ and $P'$ if necessary, we may assume $y''(-y)>0$.  Applying the first case to $P''$ and $\bar P$,
 \[
 \theta(P')=\theta(P''/P)=\theta(P'')+\theta(\bar P)=\theta(PP')-\theta(P),
 \]
 which is \eqref{eq:additivity}. This establishes existence.

\subsection{Properties of the function \texorpdfstring{\large $\theta$}{theta}}\label{subsec:theta}
{\em The function $\theta(P)$ is a nondecreasing function of $y$, for $P$ in the right-half unit circle.}

To derive this, let $P$ be in  the unit circle first quadrant. From \eqref{eq:mean},  $y_1$ is an increasing function of $y$. Similarly, with  $y_2$ playing the role of $y_1$, $y_2$ is an increasing function of $y_1$, hence an increasing function of $y$. Continuing in this manner, $y_n$, $n\ge1$, are increasing functions of $y$. By \eqref{eq:2^ny}, $\theta_n$, $n\ge1$, are increasing functions of $y$. Passing to the limit, and since 
$$\theta(P) = \theta\left(\sqrt{1-y^2},y\right)$$
is an odd function of $y$ on the right-half unit circle, the result follows.

\subsection{The number \texorpdfstring{\large $\pi$}{pi}}
Define $\pi := 2\theta(J)$ where $J$ is the point $(0,1)$. Since $2\theta_1(J)=2\sqrt{2}$ and $2\theta^*_1(J)=4$,  $2\sqrt{2}<\pi<4$. Archimedes achieved \eqref{eq:archest} by doubling a hexagon four times,  resulting in  a 96-gon. To achieve \eqref{eq:archest} using \eqref{eq:archseq}, one needs to calculate
\begin{equation*}
2\theta_6 = 64 \sqrt{2 - \sqrt{2 + \sqrt{2 + \sqrt{2 + \sqrt{2 + \sqrt{2}}}}}},
\end{equation*}
together with its circumscribed counterpart $2 \theta_6^*=2\theta_6/x_6$.
This corresponds to doubling an octagon four times, resulting in a 128-gon. 

Plugging $x=x_{n-1},y=y_{n-1}$ into \eqref{eq:mean} we obtain 
\[
(x_{n},y_{n})=
\left(\dfrac{1+x_{n-1}}{\sqrt{\mathstrut 2+2x_{n-1}}},\dfrac{y_{n-1}}{\sqrt{\mathstrut 2+2x_{n-1}}}\right) \,,
\]
from which follows $2x_{n}=\sqrt{\mathstrut 2+2x_{n-1}}$ and $2y_{n}=y_{n-1}/x_{n}$.  By \eqref{eq:2^ny}, this implies
$\theta_{n} = \theta_{n-1}/x_{n}$, which telescopes as follows at $P=J$:
\[
\theta_{n}= \dfrac{\theta_1}{x_2 \cdots x_{n}}
            = \dfrac1{x_1 \cdots x_{n}},
\]
because $\theta_1=\sqrt{2}$ and $x_1=1/\sqrt{2}$. Since $\theta_n \to \pi/2$,
\[
\dfrac2{\pi} = x_1 x_2 x_3 \cdots =\frac{\sqrt{2}}{2}
\frac{\sqrt{2+\sqrt{2}}}{2}
\frac{\sqrt{2+\sqrt{2+\sqrt{2}}}}{2} \cdots \,.
\]
This is Vi\`ete's product formula (1593), said to be ``the oldest infinite product" \cite[p. 136]{Osler-2007}
 and ``the first time an infinite process was explicitly written as a mathematical formula" \cite[p. 50]{Maor-1998}.

By \eqref{eq:archpost}, and since $\theta(P)$ is an odd function of $y$, $\theta(P)$ is continuous at $P=N$ and $\theta(P)\ne0$ when $P \ne N$. If $P\ne P'$ are in the right-half unit circle and close to each other, then $P/P' \ne N$ and close to $N$. By \eqref{eq:additivity},
$$\theta(P/P')=\theta(P)-\theta(P'),\qquad x>0,x'>0.$$ 
It follows that $\theta(P)$ is continuous and injective on the right-half unit circle. By the doubling identity \eqref{eq:trigdouble}, $\theta(P)$ is continuous and injective on the punctured unit circle.
  
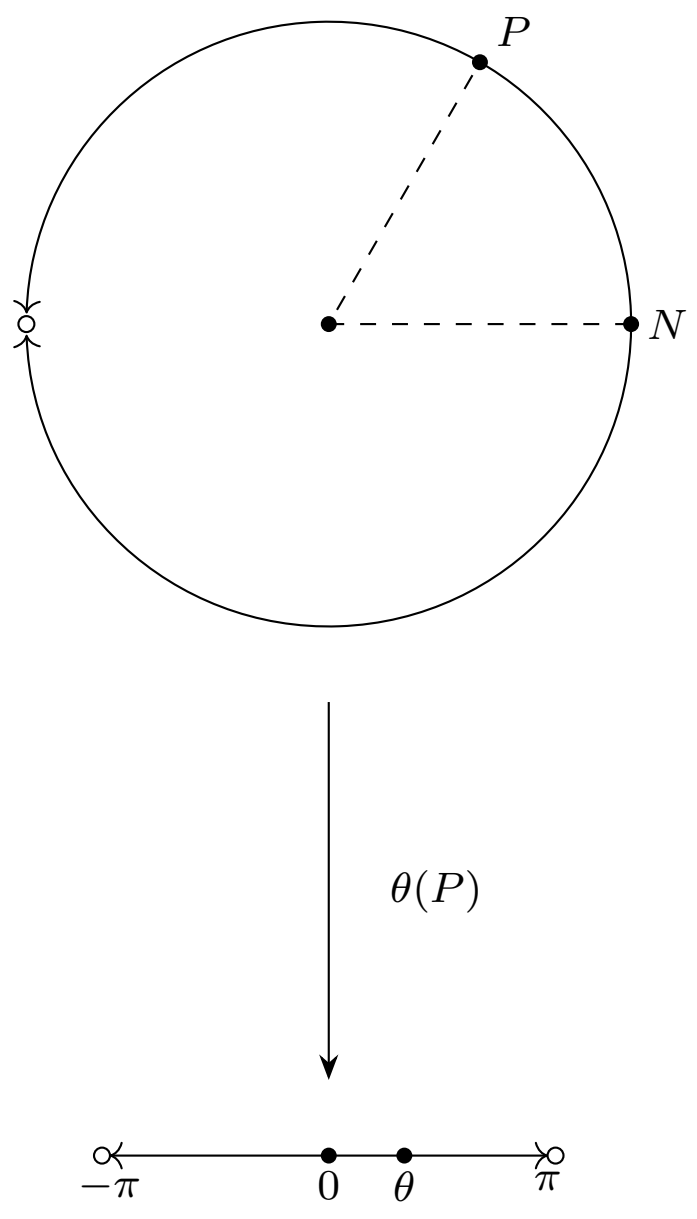
\begin{figure}[H]
	\centering
	\footnotesize
	\begin{tikzpicture}[scale=1.0]
		\coordinate (O) at (0,0);
		\coordinate (N) at (2,0);
		\coordinate (J) at (-2,0);
		\coordinate (P) at (1,1.73204);
		\draw[<->] (-178:2) arc (-178:178:2);
		\draw (J) circle [radius=1.5pt];
		\fill (P) circle [radius=1.5pt];
		\fill (N) circle [radius=1.5pt];
		\draw (P)node[above right]{$P$};
		\draw (2,0)node[right]{$N$};
		\draw[dashed] (N) -- (O) -- (P);
		\fill (O) circle [radius=1.5pt];
		\draw[-Stealth] (-0,-2.5) --node[right]{$\quad \theta(P)$} (0,-5.00);
		\begin{scope}[shift={(0,-5.5)}]
			\draw[<->] (-1.45,0)node[below]{$-\pi$} -- (1.45,0)node[below]{$\pi$};
			\draw (.5,0)node[below]{$\theta$};
			\draw (0,0)node[below]{$0$};
			\fill (0,0) circle [radius=1.5pt];
			\fill (.5,0) circle [radius=1.5pt];
			\draw[fill=white] (-1.5,0) circle [radius=1.5pt];
			\draw[fill=white] (1.5,0) circle [radius=1.5pt];
		\end{scope}
	\end{tikzpicture}
	\caption{Angle measure.}
	\label{fig:arch}
\end{figure}

 Since $\theta(P)$ is an increasing function of $y$ and $\theta(\pm J)=\pm\pi/2$, $\theta(P)$ maps  the right-half unit circle into $(-\pi/2,\pi/2)$. By the intermediate value theorem, $\theta(P)$ is a continuous bijection from the right-half  unit circle onto $(-\pi/2,\pi/2)$. By \eqref{eq:trigdouble} again,  $\theta(P)$ is a continuous bijection from the punctured unit circle $P\ne-N$ onto $(-\pi,\pi)$ (Figure \ref{fig:arch}). 

\section{Trigonometry}
\label{sec:trig}

\subsection{The fundamental theorem of trigonometry}
The basic result of trigonometry is the existence and uniqueness, up to scaling, of a continuous map  $P(\theta)=(\cos\theta,\sin\theta)$ of the real line into the unit circle satisfying the addition formula
\begin{equation}
\label{eq:trighomo}
P(\theta)P(\theta')=P(\theta + \theta')
\end{equation}
for all $\theta$, $\theta'$. This form of the addition formula is implicit in the work of de~Moivre. Ancient addition formulas are found in the Almagest~\cite{Too2}, albeit only for chord lengths.
Separating the addition formula into real and imaginary parts yields the usual trigonometric identities for $\sin\theta$ and $\cos\theta$. We already attained a finite-resolution version of this in \eqref{eq:half-chord_addition}.

By \eqref{eq:trighomo}, any such map satisfies  $P(0)=N$ and $P(-\theta)=\bar P(\theta)$. 
An immediate consequence of \eqref{eq:trighomo} is the
\emph{doubling property} 
\begin{equation}
\label{eq:trigdoubleinv}
P(\theta)=P(\theta/2)^2.
\end{equation}
If $P(\theta)$ satisfies \eqref{eq:trighomo}, then $P(\tau\theta)$ satisfies \eqref{eq:trighomo}, for any real $\tau$. Hence, at best,  there is uniqueness up to rescaling. To show this is indeed the case, we need to attach a scale to each such map $P(\theta)$. This scale, or yardstick, is its least positive period.

Given a continuous map $P(\theta)$ satisfying \eqref{eq:trighomo} for all $\theta$, $\theta'$, we say a real $\omega$ is a \emph{period}\footnote{A more conventional definition of \emph{period} requires $P(\theta+\omega)=P(\theta)$, for every $\theta$. Because of the addition formula \eqref{eq:trighomo}, $P(\theta+\omega)=P(\theta)P(\omega)$, so the two definitions are equivalent in our context.}  if $P(\omega)=N$. Then $0$ is a period, every integer multiple of a period is a period, and a limit of periods is a period.

\step{existence}  $\theta(P)$ is a continuous bijection from the punctured circle to $(-\pi,\pi)$, with inverse $P(\theta)$. Restricting $\theta(P)$ to the right half-circle and using the $y$ coordinate as a parameter, $\theta(P)$ is nondecreasing (\S\ref{subsec:theta}), hence $\theta(P)$ can be viewed as a continuous, monotone function onto $(-\pi/2,\pi/2)$; injectivity shows the monotonicity is strict. The inverse of a continuous strictly monotone function is continuous and strictly monotone. Thus $P(\theta)$ is continuous on $(-\pi/2,\pi/2)$. Using the doubling property, $P(\theta)$ is continuous on all of $(-\pi,\pi)$.

Since $\theta(J)=\pi/2$,  $P(\pi/2)=J$. By \eqref{eq:additivity},  $P(\theta)$  satisfies \eqref{eq:trighomo} on $(-\pi/2,\pi/2)$, hence \eqref{eq:trigdoubleinv} holds for $\theta$ in $(-\pi,\pi)$.
 
Fix $\alpha\ge\pi$, and suppose  $P(\theta)$ is defined on $(-\alpha,\alpha)$ and satisfies  \eqref{eq:trighomo} on $(-\alpha/2,\alpha/2)$. If $\hat P(\theta)$ extends $P(\theta)$ to $(-2\alpha,2\alpha)$ and satisfies  \eqref{eq:trighomo} on $(-\alpha,\alpha)$, then  $\hat P(\theta)=\hat P(\theta/2)^2=P(\theta/2)^2$ on $(-2\alpha,2\alpha)$, hence $\hat P(\theta)$ is uniquely determined. Conversely,  define $\hat P(\theta) :=P (\theta/2)^2$ on $(-2\alpha,2\alpha)$. Then $\hat P(\theta)$  satisfies  \eqref{eq:trighomo} on $(-\alpha,\alpha)$, since 
\[
\begin{split}
\hat P(\theta)\hat P(\theta')
&=P(\theta/2)^2P(\theta'/2)^2 \\
&=P((\theta+\theta')/2)^2
=\hat P(\theta+\theta') ,
\end{split}
\]
for $\theta$, $\theta'$ in $(-\alpha,\alpha)$, and $\hat P(\theta)$ extends $P(\theta)$, since $\hat P(\theta)=P(\theta/2)^2=P(\theta)$ on $(-\alpha,\alpha)$.

Iterating this, $P(\theta)$ can be extended uniquely  to $(-2\pi,2\pi)$, $(-4\pi,4\pi)$, $(-8\pi,8\pi)$, $\dots$, with the extensions satisfying \eqref{eq:trighomo} on $(-\pi,\pi)$, $(-2\pi,2\pi)$, $(-4\pi,4\pi)$, $\dots$, resulting in a continuous map $P(\theta)$ satisfying \eqref{eq:trighomo} for all reals. Since $P(\pi)=P(\pi/2)^2=J^2=-N$, $P(\theta)$ is surjective. Since $\theta(P)$ is bijective and $P(2\pi)=P(\pi)^2=N$, $2\pi$ is the least positive period of $P(\theta)$. We write $P(\theta)=(\cos\theta,\sin\theta)$ for this extension.
 
\step{uniqueness, step 1} Let $P(\theta)$ be any non-constant continuous map $P(\theta)$ of the real line into the unit circle satisfying \eqref{eq:trighomo}. \emph{Then $P(\theta)$ has a positive period:}

Write  $P(\theta)=(x(\theta),y(\theta))$. Since $0$ is a period, $x(0)=1$. Since $P(\theta)$ is non-constant, there is an $\alpha\ne0$ with $x(\alpha)<1$. Since $\cos0=1$, by continuity, there is an $n\ge1$ with $x(\alpha)<\cos(2\pi/n)<1$. By the intermediate value theorem, there is a $\beta\ne0$ with $x(\beta)=\cos(2\pi/n)$, hence one of $y(\pm\beta)=\pm y(\beta)$ equals $\sin(2\pi/n)$. By \eqref{eq:trighomo}, since $x(\theta)$ is even,
\[
\begin{split}
P(\pm n\beta)=P(\pm\beta)^n &= (x(\pm \beta),y(\pm\beta))^n \\
&=(\cos(2\pi/n),\sin(2\pi/n))^n \\
&=(\cos(2\pi),\sin(2\pi))=N,
\end{split}
\] 
hence one of $\pm n\beta$ is a positive period. 

\step{uniqueness, step 2}
\emph{$P(\theta)$ has a least positive period $\omega$:}

Since a limit of periods is a period, the infimum $\omega$ of all positive periods is a period. For every real $\theta$ and every period $\beta>0$, there is some integer $n$ for which  $n\beta\le\theta<(n+1)\beta$. Hence for every period $\beta>0$, every real $\theta$ lies within distance $\beta$ of some period. If $\omega=0$, there would be arbitrarily small positive periods $\beta$. This would imply every real $\theta$ is the limit of some sequence of periods, and thus $P(\theta)\equiv N$. Since identically constant functions are disallowed,  $\omega>0$. 

\step{uniqueness,  step 3}
By rescaling $P(\theta)$ to  $P(\omega\theta/2\pi)$, we may  assume $\omega=2\pi$. Then $P(\pi)=-N$, so $P(\pi/2)=J$ or $P(-\pi/2)=J$. By rescaling $P(\theta)$ to $P(-\theta)$ if necessary, we may assume  $P(\pi/2) = J$.  \emph{Then $P(\theta)$ is uniquely determined:}

Since $2\pi$ is the least positive period,  $P(\theta)\ne\pm N$ for $0<|\theta|<\pi$.   In particular, $\sqrt{P(\theta)}$ is defined for $-\pi<\theta<\pi$.   Since \eqref{eq:trighomo} implies \eqref{eq:trigdoubleinv}, 
\begin{equation}
\label{eq:sqrt2}
P(\theta/2)=\sqrt{P(\theta)},\qquad -\pi<\theta<\pi,
\end{equation}
up to sign (beneficial ambiguity again). We claim \eqref{eq:sqrt2} is correct as written, with no sign adjustment.  This is immediate when $\theta=0$, so assume $\theta\ne0$. If the imaginary parts of $P(\theta)$ and $P(\theta/2)$ have opposite signs,  then, by the intermediate value theorem, for some $\theta'$ between $\theta$ and $\theta/2$, we have $P(\theta')=N$ or $P(\theta')=-N$. Since this can't happen, the imaginary parts of $P(\theta)$ and $P(\theta/2)$ must have the same sign. It follows that  the imaginary parts of $\sqrt{P(\theta)}$ and $P(\theta/2)$ have the same sign, establishing the claim.  Using \eqref{eq:trighomo} and \eqref{eq:sqrt2} repeatedly, $P(\theta)$ satisfies 
$$
P(\theta\pi/2) = J^\theta,\qquad J=(0,1),
$$
for all dyadic rationals $\theta=k/2^n$. Since the dyadic rationals are dense and $P(\theta)$ is continuous, this determines $P(\theta)$, completing the proof of the fundamental theorem.

\subsection{Trigonometric derivatives}
Inserting $P=(\cos\theta,\sin\theta)$  in the Elbow Bounds \eqref{eq:archpost} leads to \eqref{eq:sintheta} as follows. The left inequality gives $\sin\theta < \theta$, or 
\[
\frac{\sin\theta}{\theta} <1.
\]
The right inequality gives $\theta < 1- \cos\theta + \sin\theta$, or 
\begin{equation*}
\begin{split}
0<
1-\frac{\sin\theta}{\theta} &<
\frac{1-\cos\theta}{\theta} \\ 
&= 
\frac{1-\cos\theta}{\theta}\cdot
\frac{1+\cos\theta} 
{1+\cos\theta} \\
&=
\frac{\sin^2\theta}{\theta (1+\cos\theta)} \\
&< 
\frac{\theta}{\sin\theta}
\cdot \frac{\sin^2\theta}{\theta (1+\cos\theta)} \\
&=\frac{\sin\theta} {1+\cos\theta},
\end{split}
\end{equation*}
and this last expression tends to zero as $\theta \rightarrow 0$. The reasoning here depends on the Elbow Bounds~\eqref{eq:archpost}, hence requires $0 < \theta < \pi/2$ and proves the right hand limit portion of \eqref{eq:sintheta}; that $\sin\theta/\theta$ is an even function verifies \eqref{eq:sintheta} in full. As a byproduct, we also obtain 
\[
\lim_{\theta \rightarrow 0} \dfrac{1-\cos\theta}{\theta} =0. 
\]
With this limit and \eqref{eq:sintheta} established, the derivatives of trigonometric functions can be computed in the standard way.

\section*{Afterwards}

S: ``Well, that was quite a journey you took us through. Would you regard the angle measure construction as a furthering of the G.D.~Birkhoff postulates?"

T: ``In a sense, yes. By analogy, an abstract definition of the real numbers is proposed, as a complete ordered field, and then Richard Dedekind comes along and builds a model with his Dedekind Cuts. Historically, the cuts came first, and then the abstract, unique description, though with precursors to the abstraction far earlier on. The Birkhoff postulates include the Protractor Axiom, and we have built a model for a protractor, following Archimedes. 

S: ``Yes, I've heard of models in mathematics.  This gives me a clearer idea of what they're about."

T: ``Speaking of mathematical models, my first college mathematics professor was Michael Morley. Only later did I learn that he was a distinguished model theorist. Interestingly, in the Real Analysis class that he taught, he had us read on our own about the Dedekind Cuts model in the first chapter of Rudin's \emph{Principles} book—Baby Rudin we called it, and started the course assuming everything therein."

S: ``The limit $\lim_{\theta\to 0} \frac{\sin \theta}{\theta}$ reminds me of my freshman year. We covered it in my calculus class and it seemed routine. But my roommate, who was taking the same class with a different professor kept coming back with stories. When the limit was introduced, he said that his professor `couldn't do it'. He knew that his professor was very smart, maybe even internationally distinguished, so getting stuck on a basic limit seemed out of character. That professor came back the next day and said he still couldn't do it. He looked it up in the book and said that the book was wrong; he tried other books and they were also wrong. I think they ended up just taking the value of the limit as an axiom.

Would you suggest including the angle measure development in a calculus course?"

T: ``Well, I wouldn't subject them to the whole story, but I would mention the sublime challenge of defining arc length and that this limit is intimately related to it. I would also show the Archimedes Convexity Postulate and its Elbow Bounds specialization. I think that students would appreciate the ancient historical span of calculus concepts, and the exposition challenges that even current textbooks face."

S: ``You had also mentioned that rotations lead to the concept of the integers. I remember being introduced to negative numbers in grade school via the concept on \emph{debt}. In hindsight, this is a rather abstract notion, not easily incorporated into arithmetic. Full rotations, on the other hand, forward and back, seem more intuitive and were surely known to ancient people long before debt. After all, they saw rotational motion in celestial bodies and in bird maneuvers."

T: ``Yes, quite so. And partial rotations are implicit in Archimedes' work, with limits of arcs hinting at the completeness of the real numbers."

The server refreshes the interlocutors' coffee once more.

T: ``You know, the ambiguity of angle stacking brings back memories. Years ago, in another university, when I had an administrative responsibility I was called to a high level meeting with an academic officer. Sliding off the agenda he offered: `What are you teaching our students? You make them use obscure concepts like \emph{absolute value}. In Engineering we have a much better construct: plus/minus, $\pm$. This remark didn't go over well at the time. Looking back, though, I see that he was anticipating \emph{beneficial ambiguity}."

S: ``Was Archimedes an engineer?"

\section*{Acknowledgements}
The authors thank Sheldon Axler and Daniel Klain  for helpful remarks. \\

The cafe photograph comes from Wikipedia under a Creative Commons Attribution-Share Alike 3.0 Unported license: \\
{\footnotesize
\url{commons.wikimedia.org/wiki/File:Cafe_Tortoni_La_mesa_que_Borges_acostumbraba_ocupar.jpg}}

The AI tools of Gemini, ChatGPT and Claude were used to improve the \texttt{tikz} diagrams and as sounding boards for writing style.

\end{multicols}

\noindent
\begin{minipage}[t]{0.48\textwidth}
	\footnotesize
	University of Massachusetts, Boston
	\par\medskip
	{\itshape Email address}: \upshape\ttfamily eric.grinberg@umb.edu
\end{minipage}\hfill
\begin{minipage}[t]{0.48\textwidth}
	\footnotesize
	Temple University, Philadelphia
	\par\medskip
	{\itshape Email address}: \upshape\ttfamily omar.hijab@temple.edu
\end{minipage}

\makeatletter
\let\addresses\@empty   
\makeatother


\begin{thebibliography}{18}

\bibitem{B} G. D. Birkhoff, A Set of Postulates for Plane Geometry, \textit{Ann. Math.} \textbf{33} (1932) 329--345.

\bibitem{B-R} M.~Beck and S.~Robins, Computing the Continuous Discretely: Integer-Point Enumeration in Polyhedra, Springer (2007).

\bibitem{CJ} R. Courant and F. John, \textit{Introduction to Calculus and Analysis, I}, Interscience Publishers (1965).

\bibitem{E} L. Ehrenpreis, personal communication.

\bibitem{GW} P. M. Gruber and J. M. Wills, \textit{Convexity and Its Applications}, Birkh{\"a}user Basel (1983).

\bibitem{Hardy-1908}
G. H. Hardy, \emph{A Course of Pure Mathematics}, 10th ed., Cambridge University Press, Cambridge, 1993.


\bibitem{Ha} R. Hartshorne, \textit{Geometry: Euclid and Beyond}, Springer (2000).

\bibitem{H} T. L. Heath, \textit{The Works of Archimedes}, Cambridge University Press (1897).

\bibitem{Hi} O. Hijab, \textit{Introduction to Calculus and Classical Analysis}, Springer (2016).

\bibitem{L} S. Lang, \textit{Basic Mathematics}, Springer (1988).

\bibitem{Landau-1956}
Edmund Landau, \textit{Differential and Integral Calculus}, AMS Chelsea Publishing, Providence, RI, 1956. Translated by Melvin Hausner and Martin Davis.

\bibitem{Lee-2013}
John M. Lee, \emph{Axiomatic geometry}, Pure and Applied Undergraduate Texts, vol.~21, American Mathematical Society, Providence, RI, 2013.

\bibitem{Maor-1998}
E. Maor, \emph{Trigonometric Delights}, Princeton University Press, 1998.


\bibitem{MB} E. Maor and R. W. Barnard, Trigonometry, \textit{Encyclopaedia Britannica} (2023),
{\scriptsize \url{www.britannica.com/science/trigonometry}}.

\bibitem{Osler-2007}
T.~J. Osler, \emph{A simple geometric method of estimating the error in using Vieta's product for $\pi$}, Int. J. Math. Educ. Sci. Technol. \textbf{38} (2007), no.~1, 136--142, DOI 10.1080/00207390601002799.



\bibitem{Richman} F. Richman, A Circular Argument, \emph{The College Math. Journal} \textbf{24}:2, 160--162 (1993).

\bibitem{Rudin} Walter Rudin, \textit{Real and Complex Analysis}, third ed., McGraw-Hill Book Co., New York, 1987.

\bibitem{SS} K. V. Sarma and Kripa Shankar Shukla, \textit{Aryabhatiya of Aryabhata}, India NSA (1976).

\bibitem{T} O. Toeplitz, \textit{The Calculus, A Genetic Approach}, University of Chicago Press (1963).

\bibitem{Too1} G. J. Toomer, The Chord Table of Hipparchus. \textit{Centaurus} \textbf{18} (1973) 6--28.

\bibitem{Too2} G. J. Toomer, \textit{Ptolemy's Almagest}, Princeton University Press (1998).

\end{thebibliography}
\end{document}